\documentclass[12pt,a4]{amsart}
\usepackage{amscd,amsmath,amssymb,amsfonts,ifthen,hyperref,tikz-cd}
\usepackage{calligra}
\usepackage{mathtools}

\usepackage{bbold}
\usepackage{color} 
\usepackage{graphicx}
\usepackage[all]{xypic}
\usepackage{quiver}    

\usepackage{fourier}
\usepackage{srcltx}

\newtheorem{claim}{Claim}[section]
\newtheorem{ex}[claim]{Example}
\newtheorem{thm}[claim]{Theorem}
\newtheorem{lem}[claim]{Lemma}
\newtheorem{cor}[claim]{Corollary}
\newtheorem{prop}[claim]{Proposition}
\theoremstyle{definition}
\newtheorem{defn}[claim]{Definition}

\newtheorem{rem}[claim]{Remark}

\newtheoremstyle{custombold}
  {\topsep}   
  {\topsep}   
  {\itshape}  
  {}          
  {\bfseries} 
  {.}         
  {.5em}      
  {}          

\theoremstyle{custombold}

\newtheorem{innercustomtheorem}{Theorem}
\newenvironment{customtheorem}[1]
  {\renewcommand\theinnercustomtheorem{#1}\innercustomtheorem}
  {\endinnercustomtheorem}

\newtheorem{innercustomcorollary}{Corollary}
\newenvironment{customcorollary}[1]
  {\renewcommand\theinnercustomcorollary{#1}\innercustomcorollary}
  {\endinnercustomcorollary}

\DeclareMathAlphabet{\mathcalligra}{T1}{calligra}{m}{n}
\DeclareFontFamily{OT1}{pzc}{}
\DeclareFontShape{OT1}{pzc}{m}{it}{<-> s * [1.10] pzcmi7t}{}
\DeclareMathAlphabet{\mathpzc}{OT1}{pzc}{m}{it}

\newcommand{\Spec}{\text{\rm Spec\hspace{.1ex}}}

\newcommand{\Hrm}{\mathrm{H}}

\usepackage{geometry}
\author[VQ Bao]{V\~o Qu\^oc Bao}
\email[VQ Bao]{vqbao@math.ac.vn}
\address{Institute of Mathematics, Vietnam Academy of Science and Technology}
\author[QK Nguyen]{Quang-Khai Nguyen}
\email[QK Nguyen]{nguyen@math.univ-lyon1.fr}
\address{Camille Jordan Institute, Claude Bernard University Lyon 1}

 \title{Abelian varieties are de Rham $K(\pi,1)$} 
\makeatletter
\@namedef{subjclassname@2020}{\textup{2020} Mathematics Subject Classification}
\makeatother

\thanks{The work of V\~o Qu\^oc Bao is supported by the Vingroup Innovation Foundation under grant number VINIF.2024.TS.003 and by the IMU Breakout Graduate Fellowship Program. Quang-Khai Nguyen would like to thank the Vietnam Institute of Mathematics and Labex MILYON for financial support during his master’s studies, and the European Union’s MSCA–Horizon Europe programme (grant agreement No. 101126554) for funding his doctoral studies.}

\begin{document}

\begin{abstract}
Motivated by the work of Esnault-Hai, one has the notion of \emph{de Rham  $K(\pi,1)$ schemes}, defined as follows.
Given a smooth proper geometrically connected scheme $X$ over a field $k$ of characteristic 0 and a base point $x \in X (k)$, one can define its differential fundamental group $\pi^{\mathrm{diff}}(X/k)$, which comes from the Tannakian duality of the category of coherent integrable connections on $X$. Using the formalism of $\delta$-functors, one can define natural morphisms between the group-scheme cohomology of  $\pi^{\mathrm{diff}}(X/k)$ and the de Rham cohomology of $X$. One says that $X$ with $x\in X(k)$ is {de Rham  $K(\pi,1)$} if such morphisms are all isomorphisms. In this article, we first prove that abelian varieties in characteristic $0$ are de Rham $K(\pi,1)$. In the second part of the article, we study the group-scheme cohomology of the abelianization of the differential fundamental group of a smooth proper geometrically connected scheme via its Albanese variety.
\end{abstract}
\maketitle

\tableofcontents

\section{Introduction}
Given a topological space $X$, the relation between the group cohomology of the topological fundamental group and sheaf cohomology was first studied by H. Cartan and J. Leray (see \cite[Chapter XVI.8]{CH56}). In some cases, these cohomologies are naturally isomorphic, and when this phenomenon happens, we say that $X$ is $K(\pi,1).$ A similar story has also been observed for other fundamental groups and cohomologies. For example, {\'etale $K(\pi,1)$-spaces} are schemes for which a comparison theorem holds between the group cohomology of the \'etale fundamental group and \'etale cohomology (see \cite{SGA4, Stix13,A15,A17,Abbes-Gros-2024}). In his ground-breaking work, \cite{S13}, P. Scholze showed that a connected affinoid rigid-analytic variety is a { pro-\'etale $K(\pi,1)$-space}. In \cite{EH06}, H. Esnault and P. H. Hai studied the comparison between the cohomology of the differential fundamental group and de Rham cohomology, which marked the beginning of the investigation of a de Rham analog of $K(\pi,1)$ spaces. In \cite{BHT25}, the authors showed that a smooth projective curve of positive genus is de Rham $K(\pi,1).$ In this work, we aim to continue the study in \textit{loc. cit.} by showing that abelian varieties are de Rham $K(\pi,1)$.

The notion of the differential fundamental group, first introduced in \cite{Del89} and developed further in \cite{EH06}, can be defined as follows. Let $X$ be a smooth and geometrically connected scheme over a field $k$ of characteristic $0$, and let $x \in X(k)$ be a rational point.
The \textit{differential fundamental group} $\pi^\mathrm{diff}(X/k)$ of $X$ with base point $x$ is an affine group scheme over $k,$ defined via Tannakian duality. Consider the category $\mathrm{MIC}^\mathrm{coh}(X)$ of $\mathcal{O}_X$-coherent sheaves equipped with 
$k$-linear integrable connections. This category is a $k$-linear abelian rigid tensor category equipped with the fiber functor at $x$
$$ x^*:   \mathrm{MIC}^\mathrm{coh}(X) \longrightarrow \mathrm{Vec}_k,\quad 
    (\mathcal{V}, \nabla) \mapsto \mathcal{V}_{x}, $$
which is faithfully exact, hence a neutral Tannakian category. By the Tannakian duality (see \cite[Theorem 2.11]{DM82}), this yields an affine group scheme $\pi^\mathrm{diff}(X/k)$ over $k$ such that the functor $x^\ast$
  induces an equivalence of tensor categories:
  \begin{align}\label{eq-tan1}
      x^\ast:\mathrm{MIC}^\mathrm{coh}(X) \overset{\cong}{\longrightarrow} 
 \mathrm{Rep}^\mathrm{f}(\pi^\mathrm{diff}(X/k)),
  \end{align}
where $\mathrm{Rep}^\mathrm{f}(\pi^\mathrm{diff}(X/k))$ denotes the category of finite-dimensional $k$-linear representations of $\pi^\mathrm{diff}(X/k).$ This equivalence extends to an equivalence between ind-categories (i.e., between the category of inductive direct limits of vector bundles with connection and the category of all
representations of $\pi^\mathrm{diff}(X/k)$):
\begin{equation}
\label{eq_tan2}
x^*:\mathrm{MIC}^\mathrm{ind}(X) \overset{ \cong}{\longrightarrow} \mathrm{Rep}(\pi^\mathrm{diff}(X/k)).
\end{equation}
Let $\mathrm{H}^\bullet_{\mathrm{dR}} (X,(\mathcal{V},\nabla))$ denote the de Rham cohomology of the scheme $X$ with coefficients in an integrable connection $(\mathcal{V},\nabla)$ (as defined in \cite{Kat70}, \cite{Esnault-Viehweg-1992}), and set $V:=x^*(\mathcal{V},\nabla)$. Let $\mathrm{H}^\bullet(\pi^\mathrm{diff}(X/k),V)$ denote the group-scheme cohomology of the aﬃne group scheme $\pi^\mathrm{diff}(X/k)$ with coeﬃcients in the $k$-linear representation
$V$ of $\pi^\mathrm{diff}(X/k)$, see \cite[I]{SGA3} or \cite{Jan03}.  The equivalence \eqref{eq_tan2} gives us the following morphism
\begin{align}\label{eq-tan3}
    \delta^i_{X/k}(V): \mathrm{H}^i(\pi^\mathrm{diff}(X/k),V) \longrightarrow \mathrm{H}^i_{\mathrm{dR}}(X/k,(\mathcal{V},\nabla))
\end{align}
for every $i\geq 0.$ In light of \cite{BHT25}, we introduce the following definition of \textit{de Rham} $K(\pi,1)$ schemes for such data given above, which can be seen as an analogue of the notion of $K(\pi,1)$ spaces.

\begin{defn}\label{Definition: de Rham K(pi,1)}
   A smooth geometrically connected scheme  $X$ together with a point $x\in X(k)$ is said to be \emph{de Rham} $K(\pi,1)$ if $\delta_{X/k}(V)$ is an isomorphism for all $i \geq 0 $ and for every finite-dimensional $\pi^\mathrm{diff}(X/k)$-module $V.$
\end{defn}

Achinger suggested in his comments \cite{Achinger-23} that abelian varieties should be de Rham $K(\pi,1)$. Our first result confirms Achinger's expectation.

\begin{customtheorem}{A}(\textit{Theorem \ref{Theorem: main theorem}})\label{Theorem: av is deRham K(pi,1)}
    Every abelian variety over a field of characteristic 0 is \textit{de Rham} $K(\pi,1)$.
\end{customtheorem}

\begin{rem}\label{Remark: independent of the base point}
For an abelian variety $X$, its differential fundamental group, up to isomorphism, does not depend on the choice of the base point, thanks to the translation. Thus, from now on, we will always consider the base point of an abelian variety to be its identity $0_X$.
\end{rem}

The main difficulty in proving Theorem \ref{Theorem: av is deRham K(pi,1)} is that one can not use the Poincar\'e duality as in the proof of \cite[Theorem 4.1]{BHT25} for curves. Instead, we employ a transcendental method. We first prove our result when $k=\mathbb{C}$ by translating the problem into one in algebraic topology and verifying that it holds in the topological setting. In the general case, we can find a finitely generated subfield $k'$ of $k$ for which our given data can be defined by $X'$ over $k'$. At this stage, we can compare the differential fundamental groups of $X,X'$ and $X'\otimes_{k'}\mathbb C$, and compare their corresponding cohomologies by using the Lyndon-Hochschild-Serre spectral sequence. The key idea is to exploit the fact that the differential fundamental groups of abelian varieties are commutative, from which we can decompose them into direct sums of unipotent parts and multiplicative type ones. To conclude, we apply the GAGA principle and use several base change properties of the group-scheme cohomology. 

\begin{rem}
    The transcendental method has also been used to study the homotopy sequence of the differential fundamental groups, see \cite{Zhang-13}.
\end{rem}

\begin{rem}
    The question of comparison of cohomology has also been studied in other contexts, see \cite{Proietto-Shiho-2016,Petrov-2021,Hain-2025}.
\end{rem}

Since de Rham cohomologies of abelian varieties are finite-dimensional, we obtain the following.
\begin{customcorollary}{B}[Corollary \ref{Cor-Eu1}]\label{Cor-B}
For an abelian variety of dimension $g$ over a field of characteristic $0,$ the group-scheme cohomology of its differential fundamental group is always finite-dimensional and vanishes in degrees greater than $2g.$ Moreover, the Euler characteristic of the group-scheme cohomology of its differential fundamental group is well-defined and equal to zero.
\end{customcorollary}

Note, however, that in general we do not know whether the group-scheme cohomology of the differential fundamental group of a smooth proper geometrically connected scheme is finite-dimensional or not.

Our second study in this article involves a relation between the group-scheme cohomology of the differential fundamental group of a smooth proper geometrically connected scheme and that of its associated Albanese variety.

Let $X$ be a smooth proper geometrically connected scheme over a field $k$ of characteristic 0 and $x \in X(k).$ Let $f: X \longrightarrow \mathrm{Alb_X}$ be the Albanese morphism that maps $x$ to the identity element $0_{\mathrm{Alb}_X}.$ Then we have the following functor:
$$f^{\ast}:\mathrm{MIC}^\mathrm{coh}(\mathrm{Alb}_X/k)\longrightarrow \mathrm{MIC}^{\mathrm{coh}}(X/k).$$
The pullback functor $f^\ast$ is compatible with fiber functors $x^\ast$ and $0_{\mathrm{Alb}_X}^\ast,$ hence it induces a morphism of differential fundamental groups
$$ f_{\ast}: \pi^\mathrm{diff}(X/k) \longrightarrow \pi^\mathrm{diff}(\mathrm{Alb}_X/k),$$
see \cite[Corollary 2.9]{DM82}. 
Since  the differential fundamental group $\pi^\mathrm{diff}(\mathrm{Alb}_X)$ is commutative by Lemma \ref{Lemma: diff fund group of av is commutative} below, the morphism $f_{\ast}$ factors as follows:
$$\xymatrix{ \pi^\mathrm{diff}(X) \ar@{->>}[d] \ar[r]& \pi^\mathrm{diff}(\mathrm{Alb}_X) \\
\pi^\mathrm{diff}(X)^\mathrm{ab} \ar[ur]& },$$
where $\pi^\mathrm{diff}(X)^\mathrm{ab}$ is the largest abelian quotient. 

We observe that there is an injection
$$ V^{\pi^\mathrm{diff}(\mathrm{Alb}_X/k)} \xhookrightarrow{} V^{\pi^\mathrm{diff}(X/k)^\mathrm{ab}}$$
for any $V \in \mathrm{Obj}(\mathrm{Rep}^f(\pi^\mathrm{diff}(\mathrm{Alb}_X/k))).$  Using universal $\delta$-functors (see \cite[page 205]{Har77}), we obtain the following  morphism for every $i\geq0$:
\begin{align}\label{eq-100}
    \eta^i(V): \mathrm{H}^i(\pi^\mathrm{diff}(\mathrm{Alb_X}/k),V)\longrightarrow \mathrm{H}^i(\pi^\mathrm{diff}(X/k)^\mathrm{ab},V), 
\end{align}
where $V$ is an object of $\mathrm{Rep}^f(\pi^{\mathrm{diff}}(X/k)^\mathrm{ab}).$ When $k$ is algebraically closed, we can study these maps in more detail.
\begin{customtheorem}{C}[Theorem \ref{Thm-C}]\label{customtheorem: comparision between X and its Albanese}
    Let $X$ be a smooth proper geometrically connected scheme over an algebraically closed field $k$ of characteristic $0$ and let $x\in X(k)$ be a base point. Let $f$ be the Albanese morphism mapping $x$ to $0_{\mathrm{Alb}_X}$, the identity. Then the morphism $$\pi^{\mathrm{diff}}(X/k)^{\mathrm{ab}}\to\pi^{\mathrm{diff}}(\mathrm{Alb}_X/k)$$
    is faithfully flat, and the maps
    $$\eta^i(V): \mathrm H^i(\pi^\mathrm{diff}(\mathrm{Alb}_X/k),V) \longrightarrow\mathrm H^i(\pi^\mathrm{diff}(X/k)^{\mathrm{ab}},V) $$
    are isomorphisms for all $i\geq0$ and all $V \in \mathrm{Obj}(\mathrm{Rep}^f(\pi^\mathrm{diff}(\mathrm{Alb}_X/k))).$
\end{customtheorem}

\begin{rem}
      The faithfully flat morphism $\pi^\mathrm{diff}(X/k)^{\mathrm{ab}}\xrightarrow{f_*}\pi^\mathrm{diff}(\mathrm{Alb}_X/k) $ is inspired by the work of \cite{Antei-2011,L12,BdS17,DAddezio-2023}.
\end{rem}

Our idea to prove Theorem \ref{customtheorem: comparision between X and its Albanese} comes from the complex picture, in which one can see that  $\pi^\mathrm{diff}(\mathrm{Alb}_X/k)$ and $\pi^\mathrm{diff}(X/k)^{\mathrm{ab}}$ have isomorphic unipotent parts. For general $k$, this observation can be proved by using the fact that the first de Rham cohomology of $X$ and $\mathrm{Alb}_X$ are isomorphic, whose proof itself is transcendental. We then conclude by using the Lyndon-Hochschild-Serre spectral sequence.

Since the group-scheme cohomologies of the differential fundamental group of an Albanese variety are finite-dimensional, we have the following result involving the Euler characteristic of $\pi^\mathrm{diff}(X/k)^\mathrm{ab}$. 
\begin{customcorollary}{D}[Corollary \ref{Cor-Eu2}]\label{Cor-D}
 Let $X$ be a smooth proper geometrically connected scheme over an algebraically closed field $k$ of characteristic $0$ and let $x\in X(k)$  be a base point. Then the group-scheme cohomology  $$\mathrm{H}^i(\pi^{\mathrm{diff}}(X/k)^\mathrm{ab},V)$$ is finite-dimensional for any $i\geq0$ and any $V\in \mathrm{Obj}(\mathrm{Rep}^f(\pi^{\mathrm{diff}}(X/k)^\mathrm{ab})).$ Moreover, the Euler characteristic of the group-scheme cohomology of $\pi^\mathrm{diff}(X/k)^\mathrm{ab}$  is well-defined and equal to zero.
\end{customcorollary}


Using Theorem \ref{customtheorem: comparision between X and its Albanese}, one can obtain similar isomorphisms between group cohomologies of the groups $\pi_1^{\mathrm{top}}(X,x)^{\mathrm{ab}}$ and  $\pi_1^{\mathrm{top}}(\mathrm{Alb}_X,0_{\mathrm{Alb}_X})$.
\begin{customcorollary}{E}\label{Cor-E}(Corollary \ref{Corollary: cohomology of pi top ab})
Let $X$ be a smooth proper connected complex variety, and $x$ be a base point. Let $f$ be the Albanese morphism sending $x$ to $0_{\mathrm{Alb}_X}$. Then $\pi_1^{\mathrm{top}}(\mathrm{Alb}_X,0_{\mathrm{Alb}_X})$ is the quotient of $\pi_1^\mathrm{top}(X,x)^{\mathrm{ab}}$ by its torsion subgroup, and $f$ induces isomorphisms
$$\mathrm{H}^i(\pi_1^{\mathrm{top}}(\mathrm{Alb}_X,0_{\mathrm{Alb}_X}),V)\cong\mathrm{H}^i(\pi_1^{\mathrm{top}}(X,x)^{\mathrm{ab}},V) $$
for all $i\geq0$ and all $V\in\mathrm{Obj}(\mathrm{Rep}^f(\pi_1^{\mathrm{top}}(\mathrm{Alb}_X,0_{\mathrm{Alb}_X})))$.
\end{customcorollary}

The article is organized as follows. In Section 2, we present two key lemmas on the group-scheme cohomology.
Section 3 is devoted to providing some properties of de Rham $K(\pi,1)$ schemes and their relation to the pro-algebraic completion and to $K(\pi,1)$ spaces (in the usual topological sense). In Section 4, we provide some properties of the differential fundamental groups that will be useful later on. Section 5 is devoted to
prove Theorem \ref{Theorem: av is deRham K(pi,1)} and Corollary \ref{Cor-B}. In Section 6, we prove Theorem \ref{customtheorem: comparision between X and its Albanese}, Corollary \ref{Cor-D}, and Corollary \ref{Cor-E}.

\subsection{Notations and Conventions}

\begin{enumerate}
    \item  $k$ is a field of characteristic 0.
    \item For a smooth connected scheme $X$ over $k$, let  $\mathrm{MIC}(X/k)$ denote the category of quasi-coherent
$\mathcal O_X$-sheaves equipped with $k$-linear integrable connections. 
    \item  We denote $\mathrm{MIC}^\mathrm{coh}(X/k)$ the full subcategory of $\mathrm{MIC}(X/k)$ consisting of coherent $\mathcal{O}_X$-sheaves equipped with $k$-linear integrable connections.
    \item We denote $\pi_1^\mathrm{top}(X,x)$ (resp. $\pi^\mathrm{diff}(X/k,x)$) the topological fundamental group (resp. the differential fundamental group). If the base point $x$ is clear from the context, we use $\pi^\mathrm{diff}(X/k)$ to simplify the notation.
    \item Given an aﬃne group scheme $G$ over $k,$ we let $G^\mathrm{ab}, G^\mathrm{uni}$ and $G^\mathrm{mult}$ stand respectively for the largest abelian, largest unipotent, and largest multiplicative type quotient of $G.$ When $k$ is algebraically closed, $G^\mathrm{mult}$ is also the largest diagonal quotient of $G$.
     \item  When comparing the group-scheme cohomology of the fundamental group with de Rham cohomology, we always assume that $V$ is the representation corresponding to the integrable connection $(\mathcal{V},\nabla)$ via the fiber functor.
    \item We always use the subscript to denote the result of base change (of schemes, morphisms, linear maps).
\end{enumerate}
\section{Preliminaries on the group-scheme cohomology}

    For a commutative affine group scheme $G$ over $k$, one can always decompose it as a direct product of its unipotent part and its multiplicative type part $ G = G^\mathrm{uni} \times G^{\mathrm{mult}}$, see for example \cite[Chapter IV, \S 3, no 1]{Demazure-Gabriel-1970}.

The following two lemmas are well known to experts. Due to the lack of a reference, we include proofs here.
\begin{lem}\label{Lemma: Hi G vs Hi Guni}
    Let $G$ be a commutative affine group scheme over a field $k$, and $G^{\mathrm{uni}}$ be its unipotent part. Then we have isomorphisms $$\mathrm{H}^i(G,V)\cong \mathrm{H}^i(G^{\mathrm{uni}},V^{G^{\mathrm{mult}}})$$
    for all $i\geq0$ and all finite-dimensional representations $V$ of $G$.
\end{lem}
\begin{proof}
    Using the Lyndon-Hochschild-Serre spectral sequence \cite[Proposition 6.6]{Jan03}, there is a spectral sequence whose $E_2$-page is
    $$ E_2^{p,q}=\mathrm{H}^{p}(G/G^{\mathrm{uni}},\mathrm{H}^{q}(G^{\mathrm{mult}},V))\Longrightarrow \mathrm{H}^{p+q}(G,V).
$$ We note that that taking the invariants of $G^{\mathrm{mult}}$ is an exact functor. Thus $\mathrm{H}^{q}(G^{\mathrm{mult}},V)=0$ for all $q>0$, and the spectral sequence is degenerate at the $E_2$-page. Therefore $\mathrm{H}^i(G,V)\cong \mathrm{H}^i(G^{\mathrm{uni}},V^{G^{\mathrm{mult}}})$ for all $i\geq0$ as wanted.
\end{proof}
\begin{lem}\label{Lemma: a lemma on LHS}
    Let $G$ and $H$ be commutative affine group schemes over a field $k$. Let $f:G\to H$ be a faithfully flat morphism such that $f$ induces an isomorphism on the unipotent parts. Then $f$ induces isomorphisms $$\mathrm{H}^i(H,V)\cong \mathrm{H}^i(G,V)$$
    for all $i\geq0$ and all finite-dimensional representations $V$ of $H$.
\end{lem}
\begin{proof}
    This follows from $\mathrm{H}^i(H,V)\cong \mathrm{H}^i(H^{\mathrm{uni}},V^{H^{\mathrm{mult}}})\cong \mathrm{H}^i(G^{\mathrm{uni}},V^{G^{\mathrm{mult}}})\cong \mathrm{H}^i(G,V).$
\end{proof}
\section{De Rham \texorpdfstring{$K(\pi,1)$}{K(pi,1)}}
The aim of this section is to review the notions of 
algebraically good groups, \texorpdfstring{$K(\pi,1)$}{K(pi,1)} spaces, de Rham \texorpdfstring{$K(\pi,1)$}{K(pi,1)} schemes and their relations.

\subsection{Algebraically good group}
We recall the notion of the \emph{proalgebraic completion} of an abstract group (see \cite{BLMM02}, \cite{Den23}).  
\begin{defn}
Let $\Gamma$ be an abstract group.  
The \emph{proalgebraic completion} of $\Gamma$ over a field $k$, denoted $\Gamma^{\mathrm{alg}}$, is an affine group scheme equipped with a morphism of groups
\[
  \iota \colon \Gamma \longrightarrow \Gamma^{\mathrm{alg}}(k)
\]
satisfying the following universal property:  
for any representation  $\theta: \Gamma \to \mathrm{GL}(V)$ on a finite-dimensional $k$- vector space $V,$  there is a unique representation  $\theta^\mathrm{alg}:\Gamma^{\mathrm{alg}} \to \mathrm{GL}(V)$ such that the diagram
\[
  \xymatrix{
     \Gamma \ar[r]^{\iota} \ar[rd]^\theta & \Gamma^{\mathrm{alg}}(k) \ar[d]^{\theta^\mathrm{alg}(k)} \\
     & \mathrm{GL}(V)
  }
\]
commutes.
\end{defn}

Let $\mathrm{Rep}^f(\Gamma)$ be the category of $k$-linear finite-dimensional representations of $\Gamma$ and $\mathrm{Rep}^f(\Gamma^\mathrm{alg})$ the category of $k$-linear finite-dimensional representations of the affine group scheme $\Gamma^{\mathrm{alg}}.$ According to the definition of pro-algebraic completion, the inclusion functor  
\begin{align}\label{eq-10}
    \mathrm{Rep}^f(\Gamma^\mathrm{alg})\longrightarrow \mathrm{Rep}^f(\Gamma)
\end{align}
 yields an equivalence of categories. 
\begin{ex}\label{Example: RH correspondence}
    The Riemann-Hilbert correspondence, established by Deligne \cite{Del70}, provides an equivalence of Tannakian categories
$$
\operatorname{Rep}^f\left(\pi_1^{\mathrm{top}}(X^{\mathrm{an}}, x)\right) \cong \mathrm{MIC}^\mathrm{coh}(X / \mathbb{C}),
$$
for any complex smooth proper variety ${X}$.  In other words, the differential fundamental group $\pi^\mathrm{diff}({X} / \mathbb{C})$ is the  pro-algebraic completion over $\mathbb C$ of the topological fundamental group $\pi_1^{\mathrm{top}}({X}^{\mathrm{an}}, {x}).$
\end{ex}
The equivalence \eqref{eq-10} extends to an equivalence between the ind-category 
$\mathrm{Rep}^\mathrm{ind}(\Gamma^\mathrm{alg})$ and the ind-category $\mathrm{Rep}^\mathrm{ind}(\Gamma^\mathrm{alg}).$ Since every $k$-linear representation of an affine group scheme is a direct limit of $k$-linear finite representations (see \cite{Wat79}), we have the following inclusion functor:
\begin{align}\label{eq-1}
     \mathrm{Rep}(\Gamma^\mathrm{alg}) \overset{\cong}{\longrightarrow} \mathrm{Rep}^\mathrm{ind}(\Gamma^\mathrm{alg}) \xhookrightarrow{} \mathrm{Rep}(\Gamma). 
\end{align}

For any representation $V$ of $\Gamma^{\mathrm{alg}}$, let $\mathrm{H}^\bullet(\Gamma^{\mathrm{alg}}, V)$ denote the group-scheme cohomology of the affine group scheme $\Gamma^{\mathrm{alg}}$ with coefficients in $V$ (see [\cite{SGA3},~I] or \cite{Jan03}).  
Let $\mathrm{H}^\bullet(\Gamma, V)$ denote the group cohomology of the abstract group $\Gamma$ with coefficients in $V$ (see \cite{Bro94}).  
The  equivalence \eqref{eq-1} induces, for all $i \ge 0$, natural morphisms
\[
  \iota^i \colon
  \mathrm{H}^i(\Gamma^{\mathrm{alg}}, V)
  \longrightarrow
  \mathrm{H}^i(\Gamma, V).
\]
In general, the maps $\iota^i$ are not isomorphisms, and  we are interested in the class of groups $\Gamma$ for which $\iota^i$
  is an isomorphism.

\begin{defn}\label{Definition: algebraically good group}
Let $\Gamma$ be an abstract group, and let $V$ be a finite-dimensional $k$-linear representation of $\Gamma^{\mathrm{alg}}$.  
Then $\Gamma$ is said to be \emph{algebraically good} (over $k$) if for every $i \ge 0$, the canonical map
\[
  \iota^i \colon
  \mathrm{H}^i(\Gamma^{\mathrm{alg}}, V)
  \xrightarrow{\;\cong\;}
  \mathrm{H}^i(\Gamma, V)
\]
is an isomorphism.
\end{defn}
\begin{rem}
We find that it is more convenient for us to define the definition of an algebraically good group in this way. The original one is defined differently, see \cite{KPT09}, which is equivalent to Definition \ref{Definition: algebraically good group}, cf. \cite[Lemma 4.11]{KPT09}. 
\end{rem}
\begin{ex}\label{Example: examplex of good groups}
The following groups are algebraically good over any field $k$ of characteristic $0.$
\begin{itemize}
    \item [(1)] Finite groups
    \item [(2)] Free groups of finite type
    \item [(3)] Abelian groups of finite type (\cite[Proposition 4.13]{KPT09}).
\end{itemize}
\end{ex}

\subsection{A review of \texorpdfstring{$K(\pi,1)$}{K(pi,1)} spaces}

Let $X$ be a path-connected, semi-locally simply connected topological space, and we fix a base point $x\in X.$  Let $\mathrm{Loc}(X)$ be a category of local systems on $X.$ There is a functor 
$$\rho^\ast: \mathrm{Rep}^f(\pi_1^{\mathrm{top}}(X,x)) \longrightarrow \mathrm{Loc}(X),$$ this functor yields an equivalence between these two categories  (see \cite[Corollary 1.4]{Del70}). Since $$(-)^{\pi_1^{\mathrm{top}}(X,x)} = \Gamma(X,\rho^\ast(-)),$$ using universal $\delta$-functors (see \cite[page 205]{Har77}) we obtain the following  morphism for every $i\geq0$:
\[
\rho^i : \mathrm{H}^i(\pi_1^{\mathrm{top}}(X,x), M) \longrightarrow \mathrm{H}^i(X, \rho^*M),
\]
where $M$ is a finite-dimensional representation of $\pi_1^{\mathrm{top}}(X,x).$ In general, the map $\rho^i$ is not an isomorphism; and we are interested in the class of topological spaces $X$ for which $\rho^i$ is an isomorphism for all $i\geq0$.
 
\begin{defn} 
     A topological space $X$ is called a \emph{$K(\pi,1)$ space} if $\rho^i$ is an isomorphism for every $i\geq0$.
\end{defn}

\begin{ex}\label{Example: complex tori are K(pi,1)}
For all $n\geq1$, the complex tori $(\mathbb S^1)^n$ are $K(\pi,1)$ spaces.
\end{ex}

\subsection{De Rham \texorpdfstring{$K(\pi,1)$}{K(pi,1)} schemes}
Let $X$ be a smooth geometrically connected scheme over a field $k$ of characteristic $0.$ Assume that $X$ has a rational point $x\in X(k).$  In this section, we will recap the notion of de Rham $K(\pi,1)$ schemes. 

We start by constructing the map \eqref{eq-tan3}. Let $V$ be an object of $ \mathrm{Rep}(\pi^\mathrm{diff}(X/k)).$ 
The invariant subspace of $V$ is defined to be ($k$ being equipped with
the trivial action):
\begin{align*}V^{\pi^\mathrm{diff}(X/k)}:=\mathrm{Hom}_{\pi^\mathrm{diff}(X/k)}(k,V).
\end{align*}
Since the fixed point functor $V\longmapsto V^{\pi^\mathrm{diff}(X/k)}$ is left exact, one can define the group-scheme cohomology $\mathrm{H}^{i}(\pi^\mathrm{diff}(X/k), V)$ as the $i$-th right derived functor of this functor.  Since $\mathrm{Rep}(\pi^\mathrm{diff}(X/k))$ has enough 
injectives, we get isomorphisms of right derived functors,  as stated in \cite[4.2]{Jan03}:
\begin{align}\label{33}
  \mathrm{H}^i(\pi^\mathrm{diff}(X/k),-) \cong \mathrm{Ext}^i_{\mathrm{Rep}(\pi^\mathrm{diff}(X/k))}(k, -).  
\end{align}

For a sheaf with an integrable connection $(\mathcal{V},\nabla)$ on $X/k$, 
the \textit{sheaf of horizontal sections} is defined as
\begin{align}\label{H0} 
\mathcal V^\nabla:= \mathop{\mathrm{Ker}}(\nabla : \mathcal{V} \to 
\Omega ^1_{X/k} \otimes \mathcal{V}).     
\end{align}
This is a sheaf of $k$-vector spaces. The \textit{$0$-th de Rham cohomology} of $\mathcal V$ is defined as 
$$\mathrm{H}_\mathrm{dR}^0(X/k,  (\mathcal V,\nabla)):=f_*(\mathcal V^\nabla).$$
The vector space $\mathrm{H}^0_\mathrm{dR}(X/k, (\mathcal{V},\nabla))$ can be identified with the hom-set of connections in $\mathrm{MIC}(X/k)$
$$\{\varphi:(\mathcal O_X,d) \longrightarrow (\mathcal M,\nabla)\}.$$

Since $\text{\rm MIC}(X/k)$ is equivalent to the category of left modules on the sheaf of differential operators $\mathcal{D}_{X/k}$, it has enough injectives (cf. \cite{Kat70}).
Thus, we can define the \textit{higher de Rham cohomologies} $\mathrm{H}^i_\mathrm{dR}(X/k, -)$ to be the derived functors of the functor 
$$\mathrm H^0_\mathrm{dR}(X/k, -):\mathrm{MIC}(X/k)\longrightarrow \mathrm{Vec}_k.$$ 

Recall that the differential fundamental group $\pi^\mathrm{diff}(X)$ is defined via the Tannakian duality \eqref{eq-tan1}. It follows that
$$ (-)^{\pi^{\mathrm{diff}(X/k)}} = \mathrm{Hom}_{\mathrm{Rep}(\pi^{\mathrm{diff}}(X/k))} (k, -) = \mathrm{Hom}_{\mathrm{MIC}(X/k)}((\mathcal{O}_X,d),-) =  \mathrm{H}^0_\mathrm{dR}(X/k,-),  $$
and by using the universal $\delta$-functor, we obtain the following morphisms:
\begin{align}\label{eq: contruction of delta}
    \delta^i_{X/k}(V): \mathrm{H}^i(\pi^\mathrm{diff}(X/k),V) \longrightarrow \mathrm{H}^i_{\mathrm{dR}}(X/k,(\mathcal{V},\nabla)) \text{ for all } i\geq0,
\end{align}
where $V$ is the $\pi^{\mathrm{diff}}(X/k)$-module associated with the connection $(\mathcal{V},\nabla)$.

\begin{rem}\label{Remark: isom of H1 diff and H1 dR}
The equivalence of Tannakian categories \eqref{eq-tan1} yields that $$\delta^1_{X/k}(V):\mathrm{H}^1(\pi^\mathrm{diff}(X/k),V) \longrightarrow \mathrm{H}^1_{\mathrm{dR}}(X/k,(\mathcal{V},\nabla))$$ is an isomorphism for all $V$. However, higher morphisms $\delta^i,i\geq2$, might fail to be isomorphisms, for example, when $X=\mathbb P^1$.     
\end{rem}

When $X$ with the base point $x\in X(K)$ satisfy that $\delta^{i}_{X/K}(-)$ is an isomorphism for all $i\geq0$, we say that $X$ is de Rham $K(\pi,1)$, see Definition \ref{Definition: de Rham K(pi,1)}.

\subsection{A relationship between algebraically good groups, \texorpdfstring{$K(\pi,1)$}{K(pi,1)} spaces, and de Rham \texorpdfstring{$K(\pi,1)$}{K(pi,1)} schemes.}

When $k=\mathbb C$, we can relate these notions mentioned above as follows.

\begin{prop}\label{Proposition: relations between K(pi,1) and de Rham K(pi,1)}
Let $X$ be a smooth proper connected scheme over $\mathbb{C},$ $X^\mathrm{an}$ the associated
analytic space, and $x\in X(\mathbb{C}).$ Consider the following three conditions.
\begin{itemize}
    \item [(1)]  $X$ is de Rham $K(\pi,1).$
    \item [(2)]   $X^\mathrm{an}$ is $K(\pi,1).$
    \item [(3)] The fundamental group $\pi_1^{\mathrm{top}}(X^\mathrm{an}
, x)$ is algebraically good.
\end{itemize}
Then $(1) \& (2)\Rightarrow(3)$ and $(2)\&(3)\Rightarrow(1).$
\end{prop}
\begin{proof}
    According to Deligne's Riemann-Hilbert correspondence \cite[Theorem 2.17]{Del70}, the functor 
    \begin{align*}
            \mathrm{MIC}^\mathrm{coh}(X^\mathrm{an}/\mathbb{C}) &\longrightarrow \mathrm{Loc}(X^\mathrm{an})\\
            (\mathcal{V}^\mathrm{an},\nabla^\mathrm{an}) & \mapsto \widetilde{V}:= ({\mathcal{V}^\mathrm{an}})^\nabla
    \end{align*}
    is an equivalence
between these two categories. Since $\mathrm{H}^0_{\mathrm{dR}}(X/\mathbb{C},-)=\mathrm{H}^0_{\mathrm{dR}}(X^\mathrm{an}/\mathbb{C},-) = \Gamma(X^\mathrm{an},-),$ using the universal $\delta$-functor, we obtain the following morphism:
$$ \psi^i: \mathrm{H}^i_{\mathrm{dR}}(X/\mathbb{C}, (\mathcal{V},\nabla)) \longrightarrow \mathrm{H}^i(X^\mathrm{an}, \widetilde{V}),$$
where $(\mathcal{V},\nabla)$ is an integrable connection on $X^\mathrm{an}.$

Let $(\mathcal{V},\nabla)$ be an object in $\mathrm{MIC}^\mathrm{coh}(X/\mathbb{C}).$ For every $i \geq 0,$ the morphism $\delta^i_{X/k}(V)$ (see \eqref{eq: contruction of delta}) fits into the following commutative diagram
$$
\xymatrix{\mathrm{H}^i(\pi^\mathrm{diff}(X/\mathbb{C}),V) \ar[rr]^{\delta^i_{X/k}(V)} \ar[d]^{\iota^i} & & \mathrm{H}^i_{\mathrm{dR}}(X/\mathbb{C},(\mathcal{V},\nabla)) \ar[d]^{\psi^i} \\
\mathrm{H}^i(\pi_1^{\mathrm{top}}(X^\mathrm{an},x),V) \ar[rr]^{\rho^i} & & \mathrm{H}^i(X^\mathrm{an}, \widetilde{V})}$$
since all the morphisms in the above diagram are constructed by the universal $\delta$-functor "formalism". Deligne's comparison theorem \cite[Theorem 6.12]{Del70} implies that the morphism $\psi^i$ is an isomorphism, and the result follows since $\iota^i$ is an isomorphism for all $i\geq0$ by Example \ref{Example: RH correspondence}.
\end{proof}
\begin{rem}
     Proposition \ref{Proposition: relations between K(pi,1) and de Rham K(pi,1)} can be seen as a de Rham version of \cite[Proposition 2.1.15]{A15}.
\end{rem}
\begin{cor}\label{Corollary: comlex av is de Rham K(pi,1)}
    Every abelian variety  over $\mathbb{C}$ is a de Rham $K(\pi,1)$ scheme.
\end{cor}
\begin{proof}
    Recall that every abelian variety over $\mathbb{C}$ is a complex torus, which is a $K(\pi,1)$ space by Example \ref{Example: complex tori are K(pi,1)}. In addition, its fundamental group is an abelian group of finite type, hence algebraically good by Example \ref{Example: examplex of good groups}. We then conclude by applying Proposition \ref{Proposition: relations between K(pi,1) and de Rham K(pi,1)}.
\end{proof}

\section{Differential fundamental group}
To prove Theorem \ref{Theorem: av is deRham K(pi,1)}, we need to reduce the problem to the case over $\mathbb{C}$. In order to do it, we need to understand how the differential fundamental groups behave under base change. 
\subsection{Differential fundamental groups and base change}

Let $X$ be a smooth geometrically connected scheme over a field $k$ of characteristic 0 and let $x \in X(k).$ Let $L$ be a field extension of $k$, and $X_L:=X\otimes_kL$. We have the following Cartesian diagram  
$$
\xymatrix{
  X_L \ar[d] \ar[r]^g & X \ar[d] \\
  \Spec(L) \ar[r] & \Spec(k).}
$$
Let $\mathcal{C}(X_L)$ be the full subcategory of $\mathrm{MIC}^\mathrm{coh}(X_L/L)$ whose objects, after being pushed forward along the map $X_L\xrightarrow{g} X,$ belongs to the category $\mathrm{Ind}(\mathrm{MIC}^\mathrm{coh}(X/k)).$ This $\mathcal{C}(X_L)$ is a Tannakian subcategory
and its Tannakian group is precisely $\pi^\mathrm{diff}(X/k)_{L}:=\pi^\mathrm{diff}(X/k)\times_{k}L$ thanks to \cite[10.38 and 10.41]{Del89}. We have a map  $\pi^\mathrm{diff}(X_L/L) \rightarrow \pi^\mathrm{diff}(X/k)_L$ induced by the inclusion of $\mathcal{C}(X_L)$ in  $\mathrm{MIC}^\mathrm{coh}(X_L/L)$.

\begin{lem}\label{Lemma: base change of diff fund group}
    The map $\pi^\mathrm{diff}(X_L/L) \rightarrow \pi^\mathrm{diff}(X/k)_L$ is faithfully flat. 
\end{lem}

\begin{proof}
    Since $\mathcal{C}(X_L)$ is a full subcategory of $\mathrm{MIC}^{\mathrm{coh}}(X/L)$ and is stable under taking subquotients, the result follows from the Tannakian criterion for faithfully flat morphisms in \cite{DM82}.
\end{proof}

It follows that the fixed points $(-)^{\pi^\mathrm{diff}(X/k)_{L}}$ and $ (-)^{\pi^\mathrm{diff}(X_L/L)}$ coincide.  The universal $\delta$-functor formalism yields that for every $\pi^\mathrm{diff}(X/k)_{L}$-module $V_{L}$, we have natural maps
$$ \pi_{L}^i(V_L): \Hrm^{i}(\pi^\mathrm{diff}(X/k)_{L}, V_{L}) \rightarrow \mathrm{H}^{i}(\pi^\mathrm{diff}(X_L/L),V_{L}), i\geq0.$$

\subsection{Differential fundamental groups of abelian varieties}\label{Section: Differential fundamental groups of abelian varieties}

The aim of this section is to show that the map $\pi^i_L$ constructed previously is an isomorphism when $X$ is an abelian variety. 

    \begin{lem}\label{Lemma: diff fund group of av is commutative}
Let $X$ be an abelian variety over a field $k$ of characteristic $0.$ Then the differential fundamental group $\pi^\mathrm{diff}(X/k)$ is commutative.
\end{lem}
\begin{proof}
Recall that we always consider the identity of an abelian variety as a base point. We will use an Eckmann-Hilton type argument, see \cite[Theorem 5.4.2]{EH62}. By the K\"{u}nneth theorem for $\pi^\mathrm{diff}(X/k)$ (\cite[Corollary 1.2]{Zhang-13}, \cite{Hai-2013}), the two projections from $X\times X$ to its factors induces an isomorphism
$$\pi^{\mathrm{diff}}(X\times X/k)\cong \pi^{\mathrm{diff}}(X/k)\times \pi^{\mathrm{diff}}(X/k).$$
Now, if we denote by $m$ the multiplication $X\times X\to X$ on $X$, then it induces a morphism 
$$\pi^{\mathrm{diff}}(X/k)\times\pi^{\mathrm{diff}}(X/k)\xrightarrow{m_*}\pi^{\mathrm{diff}}(X/k),$$ thus $\pi^{\mathrm{diff}}(X/k)$ is endowed with a structure of a commutative affine group scheme as wanted.
\end{proof}

\begin{prop}\label{Propostion: pi_i is isomorphism}
 Let $X$ be an abelian variety  over a field $k$ of characteristic $0.$ Let $W$ be an object in $\mathrm{Rep}_f(\pi^\mathrm{diff}(X/k)_{L}).$
Then the map 
 $$ \pi^i_{L}(W):  \mathrm{H}^i(\pi^\mathrm{diff}(X/k)_{L},W) \rightarrow \mathrm{H}^i(\pi^\mathrm{diff}(X_L/L),W)  $$
 is an isomorphism for all $i\geq 0.$
\end{prop}
    \begin{proof}
        Since $X$ is an abelian variety, the groups $\pi^\mathrm{diff}(X_L/L)$ and $\pi^\mathrm{diff}(X/k)_{L}$ are commutative by Lemma \ref{Lemma: diff fund group of av is commutative}. Thus, we have decompositions $$\pi^{\mathrm{diff}}(X_L/L)=\pi^{\mathrm{uni}}(X_L/L)\times \pi^{\mathrm{mult}}(X_L/L)\text{ and }\pi^\mathrm{diff}(X/k)=\pi^\mathrm{uni}(X/k)\times \pi^\mathrm{mult}(X/k)$$ into unipotent parts and multiplicative type ones. The faithfully flat map $ \pi^{\mathrm{diff}}(X_L/L) \longrightarrow \pi^\mathrm{diff}(X/k)_{L}$ in Lemma \ref{Lemma: base change of diff fund group} induces the faithfully flat map
$$ \pi^{\mathrm{uni}}(X_L/L) \longrightarrow \pi^\mathrm{uni}(X/k)_{L}.$$ Since $\mathbb{G}_{a}$ has no non-trivial forms, cf. \cite[14.58]{Milne-2017}, it follows from \cite[(1.2) \& (1.16)]{LM82} that $\pi^\mathrm{uni}(X_L/L )=
 (\mathbb{G}_{a})^{2g}$ where $g$ is the dimension of $X$. Since $\mathrm{H}^1(\pi^\mathrm{uni}(X/k)_{L}, L) = \mathrm{H}^1 (\pi^\mathrm{uni}(X_L/L),  L), $ it follows that $$\pi^\mathrm{uni}(X/k)_{L} =\pi^\mathrm{uni}(X_L/L) = (\mathbb{G}_a)^{2g}.$$
Using Lemma \ref{Lemma: a lemma on LHS}, we conclude that the morphism $\pi^i_{L}$ is an isomorphism for every $i\geq 0.$
    \end{proof}

\section{Abelian varieties are de Rham \texorpdfstring{$K(\pi,1)$}{K(pi,1)}}
In this section, we prove Theorem \ref{Theorem: av is deRham K(pi,1)}. Recall that we want to prove that every abelian variety (equipped with its identity element, see Remark \ref{Remark: independent of the base point}) is de Rham $K(\pi,1)$. We want to reduce the general case to the complex case by first finding a model of $X,x,(\mathcal{V},\nabla)$ over some field $k'$ embedded in both $k$ and $\mathbb{C}$, then we use the structure of the differential fundamental groups of abelian varieties to compare their group-scheme cohomology. We then can conclude by invoking the complex case Corollary \ref{Corollary: comlex av is de Rham K(pi,1)}.

\subsection{A construction}\label{Section: a general construction}
We start with a construction that aims to find our desired model over some finite generated field $k'$ over $\mathbb{Q}$.

 Let $X$ be an abelian variety over a field $k$ of characteristic $0$ equipped with a point $x\in X(k)$. There is a finite generated field $k'$ over $\mathbb{Q}$ on which $X$ and the morphisms $f, x$ are defined. In other words, we have an abelian variety $X'/k'$ that fits into the following Cartesian diagram
$$\xymatrix{
  X = X' \otimes_{k'} k \ar[d] \ar[r] & X' \ar[d] \\
  \Spec(k) \ar[r] & \Spec(k').}$$
Next, we fix an embedding $k' \xhookrightarrow{} \mathbb{C},$ so we have the following a cartesian diagram
$$\xymatrix{
  X'' = X' \otimes_{k'} \mathbb{C} \ar[d] \ar[r] & X' \ar[d] \\
  \Spec(\mathbb{C}) \ar[r] & \Spec(k').}$$
  Here, $X''$ is an abelian variety over $\mathbb{C}$.
\subsection{Proof of Theorem \ref{Theorem: av is deRham K(pi,1)}}
 We keep the same notations as in Section \ref{Section: Differential fundamental groups of abelian varieties} and Section \ref{Section: a general construction}.

\begin{cor}\label{Corollary: two base changes of diff fund groups}
    The maps $$\pi^\mathrm{diff}(X/k) \rightarrow \pi^\mathrm{diff}(X'/k')_k$$ and $$\pi^\mathrm{diff}(X''/\mathbb C) \rightarrow \pi^\mathrm{diff}(X'/k')_\mathbb{C}$$ are faithfully flat. The maps
$$ \pi_{\mathbb{C}}^i(-): \Hrm^{i}(\pi^\mathrm{diff}(X'/k')_{\mathbb{C}},-) \rightarrow \mathrm{H}^{i}(\pi^\mathrm{diff}(X''/\mathbb{C}),-)$$
and 
     $$ \pi^i_{k}(-):  \mathrm{H}^i(\pi^\mathrm{diff}(X'/k')_{k},-) \rightarrow \mathrm{H}^i(\pi^\mathrm{diff}(X/k),-)  $$
are also isomorphisms for all $i\geq0$. 
\end{cor}
\begin{proof}
    This follows directly from Lemma \ref{Lemma: base change of diff fund group} and Proposition \ref{Propostion: pi_i is isomorphism}. 
\end{proof}
\begin{lem}\label{Lemma: comparison cohomology over C and k'}
The following map
 $$\delta^i_{X'/k'}(V): \mathrm{H}^i(\pi^\mathrm{diff}(X'/k'), V)  \longrightarrow \mathrm{H}^i_{\mathrm{dR}}(X'/k', (\mathcal{V},\nabla))  $$
 is an isomorphism for all $i\geq 0$ and $(\mathcal{V},\nabla) \in \mathrm{Obj}(\mathrm{MIC}^\mathrm{coh}(X'/k')).$
\end{lem}
 \begin{proof}
Applying the flat base change theorem for de Rham cohomology (\cite[23.5]{ABC20}) and for group-scheme cohomology (\cite[Proposition 4.13]{Jan03}), we obtain the following commutative diagram
$$\xymatrix{
{\mathrm{H}^i(\pi^\mathrm{diff}(X'/k')_{\mathbb{C}},V_{\mathbb{C}})}\ar[rr]^{\cong} \ar[d]_{\pi^i_{\mathbb{C}}(V_\mathbb C)} && \mathrm{H}^i(\pi^\mathrm{diff}(X'/k'),V)\otimes \mathbb{C}\ar[dd]^{\delta^i_{X'/k'}(V)_{\mathbb C}} \\
	\mathrm{H}^i(\pi^\mathrm{diff}(X''/\mathbb{C}),V_{\mathbb{C}})\ar[d]_{\delta^i_{X''/\mathbb{C}}(V_\mathbb C)}	&&  \\
	\mathrm{H}^i_\mathrm{dR}(X''/\mathbb{C}, (\mathcal{V}_{\mathbb{C}},\nabla_{\mathbb{C}}))\ar[rr]_\cong && {\mathrm{H}^i_\mathrm{dR}(X'/k',(\mathcal{V}, \nabla)) \otimes \mathbb{C}}.
}$$
Here, $V_\mathbb{C}$ (resp. $\delta^i_{X'/k'}(V)_\mathbb C$) is the base change to $\mathbb{C}$ of $V$(resp. $\delta^i_{X'/k'}(V)$). Corollary \ref{Corollary: comlex av is de Rham K(pi,1)} yields that $\delta^i_{X''/\mathbb C}(V_\mathbb{C})$ is an isomorphism. It follows that $\delta^i_{X'/k'}(V)_\mathbb C$ is an isomorphism. Thus $\delta^i_{X'/k'}(V)$ is also an isomorphism as wanted.
 \end{proof}
\begin{thm}\label{Theorem: main theorem}
    Let $X$ be an abelian variety  over a field $k$ of characteristic $0.$ Let $(\mathcal{V},\nabla)$ be an object in $\mathrm{MIC}^\mathrm{coh}(X/k).$ Set $V:=0_X^{\ast}(\mathcal{V},\nabla).$
    Then
    \begin{align}
 \delta^{i}_{X/k}(V): \mathrm{H}^{i}(\pi^\mathrm{diff}(X/k),V) \longrightarrow \mathrm{H}_\mathrm{dR}^i(X/k, (\mathcal{V}, \nabla))
\end{align}
is an isomorphism for every $i\geq 0,$ that is, every abelian variety is de Rham $K(\pi,1)$. 
\end{thm}
 \begin{proof}
    Let $(\mathcal{V},\nabla)$ be an object of $\mathrm{MIC}^\mathrm{coh}(X/k).$ There is a finite generated field $k'$ over $\mathbb{Q}$ on which $X, (\mathcal{V},\nabla)$ and the morphisms $f, x$ are defined. In other words, there exists an integrable connection $(\mathcal{V}',\nabla')$ on $X'/k'$ such that $(\mathcal{V}',\nabla')$ is the pull back of $ (\mathcal{V},\nabla)$ via the morphism $X'\to X$. We set $V':=0_X^{\ast}(\mathcal{V}',\nabla')$, then $V=V'\otimes_{k'}k$.  Again, applying the flat base change theorem for de Rham cohomology (\cite[23.5]{ABC20}) and for group cohomology (\cite[Proposition 4.13]{Jan03}), we obtain the following commutative diagram
    $$ \xymatrix{\mathrm{H}^i(\pi^\mathrm{diff}(X/k),V) \ar[rr]^{\delta^i_{X/k}(V)} & & \mathrm{H}^i_{\mathrm{dR}}(X/k, (\mathcal{V},\nabla)) \\
  \mathrm{H}^i(\pi^\mathrm{diff}(X'/k'), V')\otimes k   \ar[u]^{\pi^i_{k}(V')} \ar[rr]^{\delta^i_{X'/k'}(V')_{k}}& &\mathrm{H}^i_{\mathrm{dR}}(X'/k', (\mathcal{V}',\nabla'))\otimes k. \ar[u]_{\cong}}$$
Corollary \ref{Corollary: two base changes of diff fund groups} yields that the map $\pi^i_{k}(V')$ is an isomorphism. Lemma \ref{Lemma: comparison cohomology over C and k'} implies that $\delta^i_{X'/k'}(V')_{k}$ is an isomorphism. We deduce that the map $\delta^i_{X/k}(V)$ is also an isomorphism as desired.
 \end{proof}
\subsection{The Euler characteristic of \texorpdfstring{$\pi^\mathrm{diff}(X/k)$}{pi-diff(X/k)}}

Let $X$ be an abelian variety of dimension $g.$ Thanks to Theorem \ref{Theorem: main theorem}, the group-scheme cohomology  $$ \mathrm{H}^i(\pi^\mathrm{diff}(X/k),V)$$
is finite-dimensional and vanishes when $i>2g.$  We may define the Euler characteristic of the group-scheme cohomology of the differential fundamental group of an abelian variety by
$$ \chi_{\pi^{\mathrm{diff}}}(V):= \sum_{i=0}^{2g}(-1)^i\dim_k\mathrm{H}^i(\pi^\mathrm{diff}(X),V),$$
where $V$ is an object in $\mathrm{Rep}^f(\pi^\mathrm{diff}(X/k)).$
\begin{cor}\label{Cor-Eu1}
    Let $X$ be an abelian variety of dimension $g$ over a field $k$ of characteristic $0.$ Let $(\mathcal{V},\nabla) \in \mathrm{Obj}(\mathrm{MIC}^\mathrm{coh}(X/k)).$ Denote $V:=0_X^{\ast}(\mathcal{V},\nabla).$ For $i\geq 0,$ then
    $$ \mathrm{H}^i(\pi^\mathrm{diff}(X/k),V)$$
is finite-dimensional and vanishes when $i>2g.$ Moreover, 
    $$\chi_{\pi^\mathrm{diff}}(V) = 0.$$
\end{cor}
\begin{proof}
It suffices to prove that $\chi_{\pi^\mathrm{diff}}(V) = 0.$ Using flat base change for de Rham cohomology (see \cite[23.5]{ABC20}), we may assume that the base field is $\mathbb{C}.$  According to \cite[Theorem 6.2, Chapter II]{Del70}, we have
$$ \mathrm{H}^i_{\mathrm{dR}}(X/\mathbb{C},(\mathcal{V},\nabla)) \cong \mathrm{H}^i(X^\mathrm{an},\tilde{V}),$$
and, combined with \cite[Proposition 2.5.4]{D04}, this implies 
$$ \chi_{\mathrm{dR}}(\mathcal{V},\nabla) = n \chi_{\mathrm{dR}}(\mathcal{O}_X,d).$$
Since $\chi_{\mathrm{dR}}(\mathcal{O}_X,d) =0$ and applying Theorem \ref{Theorem: main theorem} once again, the result follows.
\end{proof}
\section{Cohomology of the abelianization of the differential fundamental group}
In this section, we study the group scheme cohomology of $\pi^\mathrm{diff}(X)^{\mathrm{ab}}$ when $k$ is algebraically closed. To do that, we consider the Albanese morphism $f: X \longrightarrow \mathrm{Alb}_X.$ The morphism $f$ induces, for each $i\geq 0,$ a map
$$\eta^i(V): \mathrm{H}^i(\pi^\mathrm{diff}(\mathrm{Alb_X}/k),V)\longrightarrow \mathrm{H}^i(\pi^\mathrm{diff}(X/k)^\mathrm{ab},V),$$
where $V$ is an object in $\mathrm{Rep}^f(\pi^\mathrm{diff}(\mathrm{Alb}_X))$ (see \eqref{eq-100}). We aim to prove Theorem \ref{customtheorem: comparision between X and its Albanese} saying that the group-scheme cohomology of the differential fundamental group of the Albanese variety coincides with the group-scheme cohomology of $\pi^\mathrm{diff}(X)^\mathrm{ab}.$

We start with a remark on the construction of the largest abelian quotient of an affine group scheme in \cite[Page 708]{dS07}.
Every affine group scheme $G$ can be written as  $\varprojlim_{\alpha}G_\alpha $ where $\mathcal{O}(G_{\alpha}) \subseteq \mathcal{O}(G)$ is of finite type, see \cite[3.3]{Wat79}. We then define $G^\mathrm{ab}$ as the limit $\varprojlim_\alpha G^\mathrm{ab}_{\alpha}$, which has the universal property: any homomorphism $G\longrightarrow A$, where $A$ is a commutative affine group scheme, factors through $G^\mathrm{ab}.$ The affine group scheme $G^{\mathrm{ab}}$ is called the \emph{abelianization} of $G$.

\subsection{An observation over \texorpdfstring{$\mathbb{C}$}{C}}
Our proof is motivated by the following observation. 

\begin{prop}\label{Proposition: comparision with Albanese over C}
     Let $X$ be a smooth connected proper variety over a field $\mathbb C$ and $x\in X(\mathbb C)$. Let $f$ be the Albanese morphism mapping the base point $x$ to $0_{\mathrm{Alb}_X}$. Then the induced morphism 
    $$\pi^\mathrm{diff}(X/\mathbb C)^{\mathrm{ab}}\xrightarrow{f_*}\pi^\mathrm{diff}(\mathrm{Alb}_X/\mathbb C) $$is faithfully flat, and its kernel is the finite constant group scheme $(\mathrm{Tor}(\pi_1^{\mathrm{top}}(X^\mathrm{an},x)^\mathrm{ab}))^{\mathrm{alg}}$ over $\mathbb{C}$. As a consequence, we have
    $$\mathrm{H}^i(\pi^\mathrm{diff} (\mathrm{Alb}_X/k),V)\cong \mathrm{H}^i(\pi^\mathrm{diff}(X/k)^\mathrm{ab},V)$$
    for all $i\geq0$ and all $V\in\mathrm{Obj}(\mathrm{Rep}^f(\pi^{\mathrm{diff}}(\mathrm{Alb}_X/k)))$.
\end{prop}
\begin{proof}
     By the Hurewicz theorem we have $$\pi_1^{\mathrm{top}}(X^\mathrm{an},x)^\mathrm{ab}=\mathrm{H}_1(X^\mathrm{an},\mathbb Z)\text{ and }\pi_1^{\mathrm{top}}(\mathrm{Alb}_{X^\mathrm{an}},0_{\mathrm{Alb}_X})\cong \mathrm{H}_1(\mathrm{Alb}_{X^\mathrm{an}},\mathbb Z).$$
     Here, $\mathrm{Alb}_{X^\mathrm{an}}:=\dfrac{\mathrm{H}^0(X^{\mathrm{an}},\Omega^1_{X^{\mathrm{an}}})^\vee}{\mathrm{H}_1(X^{\mathrm{an}},\mathbb{Z})_{\mathrm{free}}}$ is the Albanese torus of $X^\mathrm{an}$, see \cite[11.11]{Birkenhake-Lange-2004}, whence its topological fundamental group $\pi_1^{\mathrm{top}}(\mathrm{Alb}_{X^\mathrm{an}},0_{\mathrm{Alb}_{X^\mathrm{an}}})$ is commutative. We note that $\mathrm{Alb}_{X^\mathrm{an}}$ is also the analytification of $\mathrm{Alb}_{X}$, see for instance \cite[Corollary 4.7.2.5]{HaohaoLiu-2024}. Now, the map $\mathrm{H}^1(\mathrm{Alb}_{X^\mathrm{an}},\mathbb Z)\xrightarrow{}\mathrm{H}^1(X,\mathbb Z)$ is an isomorphism. By the universal coefficient theorem, we have $$\mathrm{H}^1(X^\mathrm{an},\mathbb Z)\cong \mathrm{Hom}(\mathrm{H}_1(X^\mathrm{an},\mathbb Z),\mathbb{ Z})$$ and $$\mathrm{H}^1(\mathrm{Alb}_{X^\mathrm{an}},\mathbb Z)\cong \mathrm{Hom}(\mathrm{H}_1(\mathrm{Alb}_{X^\mathrm{an}},\mathbb Z),\mathbb{ Z}).$$ Since $\mathrm{Alb}_{X^\mathrm{an}}$ is a complex torus, $\mathrm{H}_1(\mathrm{Alb}_{X^\mathrm{an}},\mathbb{Z})$ has no torsion. We deduce that $\mathrm{H}_1({X^\mathrm{an}},\mathbb Z)\xrightarrow{}\mathrm{H}_1(\mathrm{Alb}_{X^\mathrm{an}},\mathbb Z)$, hence $\pi_1^{\mathrm{top}}(X^\mathrm{an},x)^\mathrm{ab}\xrightarrow{}\pi_1^{\mathrm{top}}(\mathrm{Alb}_{X^\mathrm{an}},0_{\mathrm{Alb}_X})$, is surjective. Moreover, its kernel is the torsion part $\mathrm{Tor}(\pi_1^{\mathrm{top}}(X^\mathrm{an},x)^\mathrm{ab})$, which is a finite commutative (abstract) group.

     Next, by Example \ref{Example: RH correspondence}, we know that $\pi^{\mathrm{diff}}(X,x)$ is the pro-algebraic completion over $\mathbb C$ of   $\pi_1^{\mathrm{top}}(X^\mathrm{an},x)$. Since taking pro-algebraic completion is right exact and commutes with taking abelianization \cite[Remark 1: (a), (c)]{BLMM02}, the map $$\pi^\mathrm{diff}(X,x)^{\mathrm{ab}}\xrightarrow{f_*}\pi^\mathrm{diff}(\mathrm{Alb}_X,0_{\mathrm{Alb}_X})$$ is also faithfully flat as desired. Next, since $\mathrm{Tor}(\mathrm{H}_1(X^{\mathrm{an}},\mathbb Z))$ is finite and commutative, it is observable subgroup\footnote{For the definition of an observable subgroup of an abstract group, we refer to \cite[Remark 1(b)]{BLMM02}.} of $\mathrm{H}_1(X^\mathrm{an},\mathbb Z)$. It follows from  \cite[Remark 1: (b)]{BLMM02} that its pro-algebraic completion $(\mathrm{Tor}(\mathrm{H}_1(X,\mathbb Z)))^{\mathrm{alg}}$ is the kernel of $$\pi^\mathrm{diff}(X,x)^{\mathrm{ab}}\xrightarrow{f_*}\pi^\mathrm{diff}(\mathrm{Alb}_X,0_{\mathrm{Alb}_X}),$$ which is a finite commutative group scheme over $\mathbb{C}$. 

     Finally, since $\pi^{\mathrm{diff}}(\mathrm{Alb}_X)$ and $\pi^{\mathrm{diff}}(X)^{\mathrm{ab}}$ are commutative, we have decompositions into unipotent parts and multiplicative type parts
$$\pi^{\mathrm{diff}}(\mathrm{Alb}_X)=\pi^{\mathrm{uni}}(\mathrm{Alb}_X)\times\pi^{\mathrm{mult}}(\mathrm{Alb}_X) \text{ and }\pi^{\mathrm{diff}}(X)^{\mathrm{ab}}=\pi^{\mathrm{diff}}(X)^{\mathrm{ab,uni}}\times \pi^{\mathrm{diff}}(X)^{\mathrm{ab,mult}}.$$
We deduce that $\pi^{\mathrm{uni}}(\mathrm{Alb}_X)\cong\pi^{\mathrm{diff}}(X)^{\mathrm{ab,uni}}$. We use Lemma \ref{Lemma: a lemma on LHS} to conclude that $\eta^i$ is an isomorphism for all $i\geq0$. 
\end{proof}

\subsection{A proof of Theorem \ref{customtheorem: comparision between X and its Albanese}}

\begin{lem}\label{Lemma: de Rham isom between X and Alb}
    Let $X$ be a smooth proper geometrically connected scheme over a field $k$ of characteristic 0 and a base point $x\in X(k)$. Let $f$ be the Albanese morphism mapping $x$ to $0_{\mathrm{Alb}_X}$, the identity. Then we have an isomorphism $$\mathrm{H}^1_{\mathrm{dR}}(\mathrm{Alb}_X,\mathcal{O}_{\mathrm{Alb}_X})\cong \mathrm{H}^1_{\mathrm{dR}}(X/k,\mathcal{O}_X).$$
\end{lem}
\begin{proof}
     Using flat base change for de Rham cohomology (see \cite[23.5]{ABC20}), we may assume that the base field is $\mathbb{C}.$ In this case, the de Rham cohomology of $X$ and $\mathrm{Alb}_X$ admit (pure) Hodge structures by \cite[Th\'eor\`eme 3.2.5]{Deligne-HodgeII}. We have the Hodge decompositions $$\mathrm{H}^1_{\mathrm{dR}}(X,\mathcal{O}_X)\cong \mathrm{H}^1(X,\mathcal{O}_X)\oplus\mathrm{H}^0(X,\Omega^1_X)$$
     and $$\mathrm{H}^1_{\mathrm{dR}}(\mathrm{Alb}_X,\mathcal{O}_{\mathrm{Alb}_X})\cong \mathrm{H}^1(X,\mathcal{O}_{\mathrm{Alb}_X})\oplus\mathrm{H}^0(X,\Omega^1_{\mathrm{Alb}_X}).$$
     By Igusa's theorem \cite{Igusa-1955}, we have an isomorphism $\mathrm{H}^0(X,\Omega^1_{\mathrm{Alb}_X})\cong \mathrm{H}^0(X,\Omega^1_X)$. By the Hodge symmetry, we also have an isomorphism $\mathrm{H}^1(X,\mathcal{O}_{\mathrm{Alb}_X})\cong \mathrm{H}^1(X,\mathcal{O}_X)$. The desired result follows.
\end{proof}

\begin{thm}\label{Thm-C}
     Let $X$ be a smooth proper geometrically connected scheme over an algebraically closed field $k$ of characteristic 0 and let $x\in X(k)$ be a base point. Let $f$ be the Albanese morphism mapping $x$ to $0_{\mathrm{Alb}_X}$. Then the morphism $$f^*:\pi^{\mathrm{diff}}(X/k)^{\mathrm{ab}}\to\pi^{\mathrm{diff}}(\mathrm{Alb}_X/k)$$
    is faithfully flat, and the maps
    $$\eta^i(V): \mathrm H^i(\pi^\mathrm{diff}(\mathrm{Alb}_X/k),V) \longrightarrow\mathrm H^i(\pi^\mathrm{diff}(X/k)^{\mathrm{ab}},V) $$
    are isomorphisms for all $i\geq0$ and all $V \in \mathrm{Obj}(\mathrm{Rep}^f(\pi^\mathrm{diff}(\mathrm{Alb}_X/k))).$
\end{thm}

\begin{proof}
First, we compare the unipotent parts.
Recall from Remark \ref{Remark: isom of H1 diff and H1 dR} that we have the isomorphisms $$\mathrm{H}^1(\pi^{\mathrm{diff}}(\mathrm{Alb}_X/k),k)\cong\mathrm{H}^1_{\mathrm{dR}}(\mathrm{Alb}_X,\mathcal{O}_{\mathrm{Alb}_X}) $$ and 
$$\mathrm{H}^1(\pi^{\mathrm{diff}}({X}/k),k)\cong\mathrm{H}^1_{\mathrm{dR}}(X/k,\mathcal{O}_X).$$
Combining with Lemma \ref{Lemma: de Rham isom between X and Alb}, we have
\begin{equation}\label{eq: isom between H1 diff Alb and H1 diff}
\mathrm{H}^1(\pi^{\mathrm{diff}}(\mathrm{Alb}_X/k),k)\cong \mathrm{H}^1(\pi^{\mathrm{diff}}({X}/k),k).    
\end{equation}
From the faithfully flat morphism $\pi^{\mathrm{diff}}(X/k)\twoheadrightarrow\pi^{\mathrm{diff}}(X/k)^{\mathrm{ab}}$, we apply the five-term exact sequence to the trivial $\pi^{\mathrm{diff}}(X/k)^{\mathrm{ab}}$-representation $k$ to obtain an injection 
\begin{equation}\label{eq: injection from H1 diff ab to H1 diff}
\mathrm{H}^1(\pi^{\mathrm{diff}}(X/k)^{\mathrm{ab}},k)\xhookrightarrow{}\mathrm{H}^1(\pi^{\mathrm{diff}}(X/k),k).    
\end{equation}
Now we would like to have an injective map $\mathrm{H}^1(\pi^{\mathrm{diff}}(\mathrm{Alb}_X/k),k)\xhookrightarrow{}\mathrm{H}^1(\pi^{\mathrm{diff}}(X/k)^{\mathrm{ab}},k)$. The group $\pi^{\mathrm{uni}}(X/k)$, which is the maximal unipotent quotient of $\pi^{\mathrm{diff}}(X/k)$, corresponds to the nilpotent Tannakian category $\mathrm{MIC}^{\mathrm{nil}}(X/k)$ of nilpotent connections via the Tannakian duality \eqref{eq-tan1}. Taking the abelianization of the faithfully flat morphism $\pi^{\mathrm{diff}}(X/k)\twoheadrightarrow\pi^{\mathrm{uni}}(X/k)$ yields a faithfully flat morphism $\pi^{\mathrm{diff}}(X/k)^{\mathrm{ab}}\twoheadrightarrow\pi^{\mathrm{uni}}(X/k)^{\mathrm{ab}}$, which in turn induces a faithfully flat morphism $\pi^{\mathrm{diff}}(X/k)^{\mathrm{ab,uni}}\twoheadrightarrow\pi^{\mathrm{uni}}(X/k)^{\mathrm{ab}}$. By \cite[Corollary 1.1.10]{Shiho-2000}, we have
\begin{equation*} 
\begin{split}
\mathrm{Lie}(\pi^{\mathrm{uni}}(X/k)^{\mathrm{ab}})^* & \cong\mathrm{Ext}^1_{\mathrm{MIC}^{\mathrm{nil}}(X/k)}((\mathcal{O}_X,d),(\mathcal{O}_X,d)) \\
 & \cong \mathrm{Ext}^1_{\mathrm{MIC}^{\mathrm{coh}}(X/k)}((\mathcal{O}_X,d),(\mathcal{O}_X,d))\\
 &\cong \mathrm{H}^1(\pi^{\mathrm{diff}}(X/k),k).
\end{split}
\end{equation*} 
and
\begin{equation*} 
\begin{split}
\mathrm{Lie}(\pi^{\mathrm{uni}}(\mathrm{Alb}_X/k))^* & \cong \mathrm{Ext}^1_{\mathrm{MIC}^{\mathrm{nil}}(\mathrm{Alb}_X/k)}((\mathcal{O}_X,d),(\mathcal{O}_X,d)) \\
 & \cong \mathrm{Ext}^1_{\mathrm{MIC}^{\mathrm{coh}}(\mathrm{Alb}_X/k)}((\mathcal{O}_X,d),(\mathcal{O}_X,d))\\
 &\cong \mathrm{H}^1(\pi^{\mathrm{diff}}(\mathrm{Alb}_X/k),k).
\end{split}
\end{equation*}
Thus we have $\mathrm{Lie}(\pi^{\mathrm{uni}}(X/k)^{\mathrm{ab}})^*\cong \mathrm{Lie}(\pi^{\mathrm{uni}}(\mathrm{Alb}_X/k))^*$. Since $\pi^{\mathrm{uni}}(X/k)^{\mathrm{ab}}$ and $\pi^{\mathrm{uni}}(\mathrm{Alb}_X/k)$ are copies of $\mathbb{G}_a$, we deduce that $\pi^{\mathrm{uni}}(X/k)^{\mathrm{ab}}\cong \pi^{\mathrm{uni}}(\mathrm{Alb}_X/k)$. Thus we have a faithfully flat morphism $\pi^{\mathrm{diff}}(X/k)^{\mathrm{ab,uni}}\twoheadrightarrow\pi^{\mathrm{uni}}(\mathrm{Alb}_X/k)$, resulting an injection
\begin{equation}\label{eq: injection from H1 diff Alb to H1 diff ab}\mathrm{H}^1(\pi^{\mathrm{diff}}(\mathrm{Alb}_X/k),k)\xhookrightarrow{}\mathrm{H}^1(\pi^{\mathrm{diff}}(X/k)^{\mathrm{ab}},k).
\end{equation}
Combining \eqref{eq: isom between H1 diff Alb and H1 diff}, \eqref{eq: injection from H1 diff ab to H1 diff}, \eqref{eq: injection from H1 diff Alb to H1 diff ab} yields that 
$\mathrm{H}^1(\pi^{\mathrm{diff}}(\mathrm{Alb}_X/k),k)\cong\mathrm{H}^1(\pi^{\mathrm{diff}}(X/k)^{\mathrm{ab}},k).$ It follows that the morphism $f^*:\pi^{\mathrm{diff}}(X/k)^{\mathrm{ab,uni}}\to \pi^{\mathrm{diff}}(\mathrm{Alb}_X/k)^{\mathrm{uni}}$ is an isomorphism.

Next we show that the morphism $f^*:\pi^{\mathrm{mult}}(X/k)\to\pi^{\mathrm{mult}}(\mathrm{Alb}_X/k)$ is faithfully flat by using \cite[Proposition 2.21]{Deligne-Milne-1982}. Via the Tannakian duality \eqref{eq-tan1}, the affine group scheme $\pi^{\mathrm{mult}}(X/k)$ corresponds to the semisimple category $\mathrm{MIC}^{{\mathrm
s}{\mathrm
s}}(X/k)$ generated by simple integrable connections on $X$ by \cite[Proposition 2.23]{Deligne-Milne-1982}. Since every simple object in $\mathrm{Rep}^f(\pi^{\mathrm{mult}}(X/k))$ is one-dimensional (recall that $k$ is algebraically closed), every simple connection on $X$ is of rank 1. We also have a similar observation for $\mathrm{Alb}_X$. Now we check the conditions in \cite[Proposition 2.21]{Deligne-Milne-1982} for  $f^*:\mathrm{MIC}^{{\mathrm
s}{\mathrm
s}}(\mathrm{Alb}_X/k)\to\mathrm{MIC}^{{\mathrm
s}{\mathrm
s}}(X/k)$.
\begin{itemize}
    \item Since the simple objects on both sides are of rank 1, $f^*$ is closed under taking
subquotients.
    \item We are left to show that $f^*$ is fully faithful. It suffices to check this condition for simple objects $(\mathcal{V,\nabla}_1)$ and $(\mathcal{W,\nabla}_2)$ in $\mathrm{MIC}^{{\mathrm
s}{\mathrm
s}}(\mathrm{Alb}_X/k)$. We observe that $$\mathrm{Hom}_{\mathrm{MIC}^{{\mathrm
s}{\mathrm
s}}(\mathrm{Alb}_X/k)}((\mathcal{V,\nabla}_1),(\mathcal{W,\nabla}_2))=k$$ if $(\mathcal{V,\nabla}_1)\cong(\mathcal{W,\nabla}_2)$, and $=0$ otherwise. Thus it reduces to checking that $f^*(\mathcal{V,\nabla}_1)\not\cong f^* (\mathcal{W,\nabla}_2)$  when $(\mathcal{V,\nabla}_1)\not\cong(\mathcal{W,\nabla}_2)$. Since every line bundle with an integrable connection on $\mathrm{Alb}_X$ belongs to $\mathrm{Pic}^0(\mathrm{Alb}_X)$ (as $\mathrm{Alb}_X$ is an abelian variety) and $\mathrm{Pic}^0(\mathrm{Alb}_X)\cong\mathrm{Pic}^0(X)$ (thanks to the universal property of $\mathrm{Alb}_X$), it remains to consider the case $\mathcal{W}=\mathcal{V}$. This case is done since the difference $\nabla_1-\nabla_2$ belongs to $\mathrm{H}^0(\mathrm{Alb}_X,\Omega^1_{\mathrm{Alb}_X})\cong\mathrm{H}^0(X,\Omega^1_{X})$.  
\end{itemize}
Therefore, the morphism  $f^*:\pi^{\mathrm{mult}}(X/k)\to\pi^{\mathrm{mult}}(\mathrm{Alb}_X/k)$ is faithfully flat. Combining with the faithfully flat morphism  $f^*:\pi^{\mathrm{uni}}(X/k)\to\pi^{\mathrm{uni}}(\mathrm{Alb}_X/k)$, we deduce that $f^*:\pi^{\mathrm{diff}}(X/k)\to\pi^{\mathrm{diff}}(\mathrm{Alb}_X/k)$ is faithfully flat. It follows that $f^*:\pi^{\mathrm{diff}}(X)^{\mathrm{ab}}\to\pi^{\mathrm{diff}}(\mathrm{Alb}_X)$ is faithfully flat, which induces an isomorphism on the unipotent parts as wanted. We conclude the proof by using Lemma \ref{Lemma: a lemma on LHS}. 
\end{proof}

\begin{rem}
    We refer readers to \cite{Brion-2018} for some properties of integrable connections on an abelian variety in arbitrary characteristic.
\end{rem}

\begin{cor}\label{Cor-Eu2}
     Let $X$ be a smooth proper geometrically connected scheme over an algebraically closed field $k$ of characteristic $0$ and let $x\in X(k)$ be a base point. Then the group-scheme cohomology  $$\mathrm{H}^i(\pi^{\mathrm{diff}}(X/k)^\mathrm{ab},V)$$ is finite-dimensional for any $i\geq0$ and any $V\in \mathrm{Obj}(\pi^{\mathrm{diff}}(X/k)^\mathrm{ab}).$ Moreover, the Euler characteristic of the group-scheme cohomology of $\pi^\mathrm{diff}(X/k)^\mathrm{ab}$  is well-defined and equal to zero.
\end{cor}
\begin{proof}
We first prove the finiteness of group-scheme cohomology of $\pi^\mathrm{diff}(X/k)^\mathrm{ab}).$ According to Theorem \ref{Thm-C}, we have 
$$ \mathrm{H}^1(\pi^\mathrm{diff}(\mathrm{Alb}_X/k),k)  \cong \mathrm{H}^1(\pi^\mathrm{diff}(X/k)^\mathrm{ab},k).$$
This implies that the unipotent part of $\pi^\mathrm{diff}(X/k)^{\mathrm{ab}}$ is $(\mathbb{G}_a)^{2g}$ where $g$ is the dimension of the Albanese variety. Thanks to Lemma \ref{Lemma: Hi G vs Hi Guni}, we obtain
$$ \mathrm{H}^i(\pi^\mathrm{diff}(X/k)^\mathrm{ab},V) \cong \mathrm{H}^i((\mathbb{G}_a)^{2g},V^{\pi^\mathrm{diff}(X/k)^\mathrm{ab,mult}}).$$
 Using Corollary \ref{Cor-B}, the desired result follows.
\end{proof}

\subsection{Cohomology of the abelianization of the topological fundamental group}
In this section, we prove Corollary \ref{Cor-E}.
\begin{cor}\label{Corollary: cohomology of pi top ab}
    Let $X$ be a smooth connected proper complex variety, and $x$ be a base point. Let $f$ be the Albanese morphism mapping $x$ to $0_{\mathrm{Alb}_X}$. Then $\pi_1^{\mathrm{top}}(\mathrm{Alb}_X,0_{\mathrm{Alb}_X})$ is the quotient of $\pi_1^\mathrm{top}(X,x)^{\mathrm{ab}}$ by its torsion subgroup, and $f$ induces isomorphisms
$$\mathrm{H}^i(\pi_1^{\mathrm{top}}(\mathrm{Alb}_X,0_{\mathrm{Alb}_X}),V)\cong\mathrm{H}^i(\pi_1^{\mathrm{top}}(X,x)^{\mathrm{ab}},V) $$
for all $i\geq0$ and all $V\in\mathrm{Obj}(\mathrm{Rep}^f(\pi_1^{\mathrm{top}}(\mathrm{Alb}_X,0_{\mathrm{Alb}_X})))$.
\end{cor}
\begin{proof}
The first statement follows from the proof of Proposition \ref{Proposition: comparision with Albanese over C}. For the second, in view of Theorem \ref{Theorem: av is deRham K(pi,1)}, Theorem \ref{customtheorem: comparision between X and its Albanese} and Proposition \ref{Proposition: relations between K(pi,1) and de Rham K(pi,1)}, it suffices to show that the groups $\pi^\mathrm{top} (\mathrm{Alb}_X,0_{\mathrm{Alb}_X})$ and $\pi^\mathrm{top}(X,x)^\mathrm{ab}$ are algebraically good. Since $X$ is a finite CW complex, the group $\pi^\mathrm{top}(X,x)$ is finitely presented. Therefore, $\pi^\mathrm{top}(X,x)^{\mathrm{ab}}$ is abelian of finite type, hence algebraically good by Example \ref{Example: examplex of good groups}. Since $\mathrm{Alb}_X$ is an abelian variety, $\pi^\mathrm{top} (\mathrm{Alb}_X,0_{\mathrm{Alb}_X})$ is also algebraically good. The corollary is proved.
\end{proof}

\begin{center}\bf 
Acknowledgments 
\end{center}
We thank Francesco Baldassarri, Michel Brion, H\'el\`ene Esnault, Phung Ho Hai, Dao Van Thinh, and Tran Phan Quoc Bao for insightful discussions. We sincerely thank Michel Brion for his careful reading of and suggestions on the first draft. The second author is grateful to Boris Adamczewski and Charles Favre for their constant support.

\bibliographystyle{alpha}
\bibliography{references}
\end{document}